\documentclass[10pt]{amsart}

\usepackage{amssymb}
\usepackage{verbatim}

\begin{document}

\newtheorem{theorem}{Theorem}    
\newtheorem{proposition}[theorem]{Proposition}
\newtheorem{conjecture}[theorem]{Conjecture}
\def\theconjecture{\unskip}
\newtheorem{corollary}[theorem]{Corollary}
\newtheorem{lemma}[theorem]{Lemma}
\newtheorem{observation}[theorem]{Observation}
\theoremstyle{definition}
\newtheorem{definition}{Definition}
\newtheorem{remark}{Remark}
\newtheorem{question}{Question}
\def\thequestion{\unskip}
\newtheorem{example}{Example}
\def\theexample{\unskip}
\newtheorem{problem}{Problem}

\numberwithin{theorem}{section}
\numberwithin{definition}{section}
\numberwithin{equation}{section}
\numberwithin{remark}{section}

\def\intprod{\mathbin{\lr54}}
\def\reals{{\mathbb R}}
\def\integers{{\mathbb Z}}
\def\complex{{\mathbb C}\/}
\def\naturals{{\mathbb N}\/}
\def\distance{\operatorname{distance}\,}
\def\rank{\operatorname{rank}\,}
\def\degree{\operatorname{degree}\,}
\def\dim{\operatorname{dim}\,}
\def\kernel{\operatorname{kernel}\,}
\def\Span{\operatorname{span}\,}
\def\codim{\operatorname{codim}}
\def\ZZ{ {\mathbb Z} }
\def\e{\varepsilon}
\def\p{\partial}
\def\rp{^{-1}}
\def\Re{\operatorname{Re\,} }
\def\Im{\operatorname{Im\,} }
\def\ov{\overline}
\def\bx{{\bf{x}}}
\def\eps{\varepsilon}
\def\ltwo{L^2}
\def\bg{{\text{big}}}
\def\smll{{\text{small}}}
\def\larrow{\overset{\rightharpoonup}{\ell}}
\def\esssup{\operatorname{ess\,\,sup}\,}

\def\scriptx{{\mathcal X}}
\def\scriptb{{\mathcal B}}
\def\scripta{{\mathcal A}}
\def\scriptk{{\mathcal K}}
\def\scriptd{{\mathcal D}}
\def\scriptp{{\mathcal P}}
\def\scriptl{{\mathcal L}}
\def\scriptv{{\mathcal V}}
\def\scripti{{\mathcal I}}
\def\scripth{{\mathcal H}}
\def\scriptm{{\mathcal M}}
\def\scripte{{\mathcal E}}
\def\scriptt{{\mathcal T}}
\def\scriptb{{\mathcal B}}
\def\frakg{{\mathfrak g}}
\def\frakG{{\mathfrak G}}
\def\boldp{{\mathbf p}}

\newcommand{\norm}[1]{ \left|  #1 \right| }
\newcommand{\Norm}[1]{ \left\|  #1 \right\| }
\newcommand{\ifof}{\Leftrightarrow}
\newcommand{\set}[1]{ \left\{ #1 \right\} }
\newcommand{\floor}[1]{\lfloor #1 \rfloor }

\author{Jonathan Bennett}
\address{School of Mathematics, University of Birmingham, Birmingham, B15 2TT, UK}
\email{J.Bennett@bham.ac.uk}

\author{Anthony Carbery}
\address{School of Mathematics, University of Edinburgh, Edinburgh, EH9 3JZ, UK}
\email{A.Carbery@ed.ac.uk}

\author{Michael Christ}
\address{
        Department of Mathematics\\
        University of California \\
        Berkeley, CA 94720-3840, USA}
\email{mchrist@math.berkeley.edu}
\thanks{The third author was supported in part by NSF grant DMS-040126}

\author{Terence Tao}
\address{Department of Mathematics, UCLA, Los Angeles CA 90095-1555, USA}
\email{tao@math.ucla.edu}

\date{May 31, 2005} 

\title
[H\"older-Brascamp-Lieb multilinear inequalities] 
{Finite bounds for H\"older-Brascamp-Lieb multilinear inequalities}

\begin{abstract}
A criterion is established for the validity of multilinear inequalities
of a class considered by Brascamp and Lieb,
generalizing well-known inequalties of H\"older, Young, and Loomis-Whitney.
\end{abstract}
\maketitle


\section{Formulation}

Consider multilinear functionals 
\begin{equation} \label{defn}
\Lambda(f_1,f_2,\cdots,f_m)
= \int_{\reals^n} \prod_{j=1}^m f_j(\ell_j(y))\,dy
\end{equation}
where 
each $\ell_j:\reals^n\to\reals^{n_j}$ is a surjective linear transformation,
and 
$f_j: \reals^{n_j}\to[0,+\infty]$. 
Let $p_1,\cdots,p_m\in[1,\infty]$.
For which $m$-tuples of exponents and linear transformations
is
\begin{equation} \label{maininequality}
\sup_{f_1,\cdots,f_m}
\frac{ \Lambda(f_1,f_2,\cdots,f_m)}
{\prod_j \|f_j\|_{L^{p_j}}}
<\infty?
\end{equation}
The supremum is taken over all $m$-tuples of nonnegative Lebesgue 
measurable functions $f_j$ having positive, finite norms.
If $n_j=n$ for every index $j$
then \eqref{maininequality} is 
essentially a restatement of  H\"older's inequality.
Other well-known particular cases include Young's inequality for convolutions
and the Loomis-Whitney inequality \cite{loomiswhitney}.

In this paper we characterize finiteness of the supremum \eqref{maininequality}
in linear algebraic terms,
and discuss certain variants and a generalization. 
In this level of generality, the question was to our knowledge
first posed by Brascamp and Lieb \cite{BL}.
A primitive version of the problem 
involving Cartesian product rather than linear algebraic structure
was posed and solved by Finner \cite{finner}; see \S\ref{section:finner} below.
In the case when the dimension $n_j$ of each target space equals
one, Barthe \cite{barthe} characterized \eqref{maininequality}. 
Carlen, Lieb and Loss \cite{CLL} gave an alternative proof for that case.
They developed an inductive analysis closely related to that of Finner,
and introduced the pivotal concept of a critical subspace.
Our analysis is a further development of those ideas.

An alternative line of analysis exists.
Although rearrangement inequalities such as that
of Brascamp, Lieb, and Luttinger \cite{bll} do not apply when
the target spaces have dimensions greater than one,
Lieb  \cite{lieb} nonetheless showed that 
the supremum in \eqref{maininequality} 
equals the supremum over all $m$-tuples of Gaussian functions,\footnote{
This situation should be contrasted with that of multilinear
{\em operators} of the same general form, mapping $\otimes_j L^{p_j}$
to $L^q$. 
When $q\ge 1$, such multilinear operators
are equivalent by duality to multilinear forms $\Lambda$.  
This is not so for $q<1$, and Gaussians are then quite far from being extremal
\cite{christtrilinear}.} 
meaning those of the form
$f_j=\exp(-Q_j(y,y))$ for some positive definite quadratic form $Q_j$.
See \cite{CLL} and references cited there for more on this approach.
In a companion paper \cite{bcct} we have given other proofs of our
characterization of \eqref{maininequality},
by using heat flow to continuously deform arbitrary functions $f_j$
to Gaussians while increasing the ratio in \eqref{maininequality}. 
This approach also extends work of Carlen, Lieb, and Loss
\cite{CLL} via a method they introduced.

\section{Results}
Denote by $\dim(V)$ the dimension of a vector space $V$.
It is convenient to reformulate the problem in a more invariant
fashion. Let $H,H_1,\dots,H_m$ be Hilbert spaces of finite, positive
dimensions.
Each is equipped with a canonical Lebesgue measure, by choosing
orthonormal bases, thus obtaining identifications with $\reals^{\dim(H)}$,
$\reals^{\dim(H_j)}$. 
Let $\ell_j:H\to H_j$ be surjective linear mappings.
Let $f_j:H_j\to\reals$ be nonnegative.
Then $\Lambda(f_1,\cdots,f_m)$ equals $\int_H \prod_{j=1}^m f_j\circ\ell_j(y)\,dy$.

\begin{theorem} \label{thm:global}
For $1\le j\le m$
let $H,H_j$ be Hilbert spaces of finite, positive dimensions.
For each index $j$ let $\ell_j:H\to H_j$ be surjective linear transformations,
and let $p_j\in[1,\infty]$.
Then \eqref{maininequality} holds 
if and only if 
\begin{equation} \label{homognecessary}
\dim(H)= \sum_j p_j\rp \dim(H_j)
\end{equation}
and
\begin{equation} \label{globalmainhypothesis}
\dim(V)
\le
\sum_j p_j\rp\dim(\ell_j(V))
\ \text{ for every subspace $V\subset H$}.
\end{equation}
\end{theorem}
\noindent
This equivalence is established by other methods in \cite{bcct}, 
Theorem~1.15.

The necessity of \eqref{homognecessary}
follows from scaling: if $f_j^\lambda(x_j) = g_j(\lambda x_j)$ 
for each $\lambda\in\reals^+$ then $\Lambda(\{f_j^\lambda\})$
is proportional to $\lambda^{-\dim(H)}$, while
$\prod_j \|f_j^\lambda\|_{p_j}$ is proportional 
to $\prod_j \lambda^{-\dim(H_j)/p_j}$.
That \eqref{globalmainhypothesis} is also necessary will be shown in 
\S\ref{gutsection}
in the course of the proof of the more general Theorem~\ref{thm:gut}.

Throughout the paper,
$\codim_W(V)$ will denote the codimension of a subspace $V\subset W$ in $W$.
Given that \eqref{homognecessary} holds,
the hypothesis \eqref{globalmainhypothesis} can be equivalently restated
as \eqref{localhypothesis}: 
$\codim_H(V)\ge \sum_j p_j\rp \codim_{H_j}(\ell_j(V))$; 
any two of these three conditions 
\eqref{homognecessary}, \eqref{globalmainhypothesis}, \eqref{localhypothesis}
imply the third.
As will be seen through the discussion of variants below,
\eqref{globalmainhypothesis} expresses a necessary condition
governing large-scale geometry (compare Theorem~\ref{thm:amalgam}), 
while \eqref{localhypothesis} expresses a necessary condition
governing small-scale geometry (compare Theorem~\ref{thm:local}). 
See also the discussion of necessary conditions for Theorem~\ref{thm:gut}.

\begin{remark}
$\Lambda$ can be alternatively expressed as a constant multiple of
$\int_\Sigma \prod_j f_j\,d\sigma$,
where $\Sigma$ is a linear subspace of $\oplus_j H_j$
and $\sigma$ is Lebesgue measure on $\Sigma$.
More exactly, $\Sigma$ is the range of the map $H\owns x\mapsto
\oplus_j \ell_j(x)$.
Denote by $\pi_j$ the restriction to $\Sigma$ of the natural projection
$\pi_j:\oplus_i H_i\to H_j$.
Then condition \eqref{globalmainhypothesis} can be restated as
\begin{equation} \label{Sigmanonsupercritical}
\dim(\tilde\Sigma)
\le
\sum_j p_j\rp \dim(\pi_j(\tilde\Sigma)) 
\text{ for every linear subspace $\tilde\Sigma\subset\Sigma$.} 
\end{equation}
\end{remark}

A local variant is also natural. Consider
\begin{equation} \label{localdef}
\Lambda_{\rm loc}(f_1,\cdots,f_m)
= \int_{\{y\in H: |y|\le 1\}}
\prod_j f_j\circ\ell_j(y)\,dy.
\end{equation} 

\begin{theorem} \label{thm:local}
Let $H,H_j,\ell_j$, and $f_j:H_j\to[0,\infty)$ be
as in Theorem~\ref{thm:global}.
A necessary and sufficient condition for there to exist $C<\infty$
such that 
\begin{equation} \label{mainlocalinequality}
\Lambda_{\rm loc}(f_1,\cdots,f_m)
\le C\prod_j\|f\|_{L^{p_j}}
\end{equation} 
for all nonnegative measurable functions $f_j$
is that for every subspace $V$ of $H$, 
\begin{equation} \label{localhypothesis}
\codim_H(V)
\ge \sum_j p_j\rp 
\codim_{H_j}(\ell_j(V)).
\end{equation} 
\end{theorem} 
\noindent This is equivalent to Theorem~8.17 of \cite{bcct},
proved there by a different method.

Certain cases of Theorem~\ref{thm:local} follow from Theorem~\ref{thm:global};
if there exist exponents $r_j$ 
satisfying the hypotheses \eqref{homognecessary} 
and \eqref{globalmainhypothesis} of Theorem~\ref{thm:global},
such that $r_j\le p_j$ for all $j$,
then the conclusion of Theorem~\ref{thm:local}
follows directly from that of Theorem~\ref{thm:global} by
H\"older's inequality, since $\Norm{f_j}_{L^{r_j}}\le C'\Norm{f_j}_{L^{p_j}}$.
But not all cases of Theorem~\ref{thm:local} are subsumed
in Theorem~\ref{thm:global} in this way. See Remark~\ref{finnerremark}
for examples.

The next theorem, in which some but not necessarily all coordinates of $y$
are constrained to a bounded set, unifies 
Theorems~\ref{thm:global} and \ref{thm:local}.
\begin{theorem} 
\label{thm:gut}
Let $H,H_0,\cdots,H_m$ be finite-dimensional Hilbert spaces
and assume that $\dim(H_j)>0$ for all $j\ge 1$.
Let $\ell_j: H\to H_j$ be linear transformations
for $0\le j\le m$, which are surjective for all $j\ge 1$.
Let $p_j\in[1,\infty]$ for $1\le j\le m$.
Then there exists $C<\infty$ such that 
\begin{equation} \label{gutinequality}
\int_{\{y\in H: |\ell_0(y)|\le 1\}}
\prod_{j=1}^m f_j\circ\ell_j(y)\,dy
\le C
\prod_{j=1}^m \|f_j\|_{L^{p_j}}
\end{equation}
for all nonnegative Lebesgue measurable functions $f_j$
if and only if 
\begin{alignat}{2}
\dim(V)&\le \sum_{j=1}^m p_j\rp \dim(\ell_j(V)) 
\ \ &&\text{for all subspaces $V\subset\kernel(\ell_0)$}
\label{firstnec}
\\
\intertext{and}
\codim_H(V)
&\ge 
\sum_{j=1}^m p_j\rp \codim_{H_j}(\ell_j(V)) 
\ \ &&\text{for all subspaces $V\subset H$.}
\label{secondnec}
\end{alignat}
\end{theorem}

This subsumes Theorem~\ref{thm:local}, by taking $H_0=H$ and
$\ell_0:H\to H$ to be the identity; \eqref{firstnec} 
then only applies to $\{0\}$, for which it holds automatically, so
that the only hypothesis is then \eqref{secondnec}.
On the other hand, Theorem~\ref{thm:global} is the special case $\ell_0\equiv 0$ 
of Theorem~\ref{thm:gut}.
In that case $\kernel(\ell_0)=H$, so 
\eqref{firstnec} becomes \eqref{globalmainhypothesis}.
In addition, the case $V=\{0\}$ of \eqref{secondnec} yields
the reverse inequality $\dim(H)\ge \sum_j p_j\rp \dim(H_j)$. 
Thus the hypotheses
of Theorem~\ref{thm:gut} imply those of Theorem~\ref{thm:global}
when $\ell_0\equiv 0$.
The converse implication also holds, 
as was pointed out in the discussion of Theorem~\ref{thm:local}.

Our next result is one of several possible discrete analogues.
\begin{theorem} \label{thm:discretegroups}
Let $G$ and $\{G_j: 1\le j\le N\}$ be finitely generated Abelian groups.
Let $\varphi_j:G\to G_j$ be homomorphisms whose ranges are subgroups
of finite indices. Let $p_j\in[1,\infty]$. Then 
there exists $C<\infty$ such that
\begin{equation}
\sum_{y\in G} \prod_{j=1}^N f_j\circ\varphi_j(y)
\le C\prod_j \|f_j\|_{\ell^{p_j}(G_j)}
\text{ for all nonnegative functions $f_j$}
\end{equation}
if and only if 
\begin{equation} \label{discretehypothesis}
\rank(H) \le \sum_j p_j\rp\rank(\varphi_j(H))
\text{ for every subgroup $H$ of $G$}.
\end{equation}
\end{theorem}
\noindent
Here the $\ell^{p_j}$ norms are of course defined with respect to
counting measure.
The constant $C$ depends of course on the torsion subgroups of the groups $G_j$. 

In $\reals^d$, for each $n\in\integers^d$
define $Q_n=\{x\in\reals^d: |x-n|\le\sqrt{d}\}$.
The space $\ell^p(L^\infty)(\reals^d)$
is the space of all $f\in L^\infty(\reals^d)$ 
for which the norm $(\sum_{n\in\integers^d} \|f\|_{L^\infty(Q_n)}^p)^{1/p}$
is finite.
\begin{theorem} \label{thm:amalgam}
Let $\ell_j:\reals^d\to\reals^{d_j}$ be surjective linear transformations.
Let $p_j\in[1,\infty]$. Then there exists $C<\infty$ such that
\begin{equation}
\int_{\reals^d} \prod_{j=1}^N f_j\circ\ell_j(y)\,dy
\le C\prod_j \|f_j\|_{\ell^{p_j}(L^\infty)(\reals^{d_j})}
\text{ for all nonnegative functions $f_j$}
\end{equation}
if and only if for every subspace $V\subset\reals^d$,
\begin{equation}
\dim(V) \le \sum_j p_j\rp \dim(\ell_j(V)). 
\end{equation}
\end{theorem}
\noindent A related result is Corollary~8.11 of \cite{bcct}.

We have assumed in all these theorems that all exponents
satisfy $p_j\ge 1$. In Theorems~\ref{thm:global}, \ref{thm:local}, and
\ref{thm:gut}, the inequalities in question are false if some
$p_j<1$. To see this, fix one index $j$. 
Take $f_i$ to be the characteristic
function of a fixed ball centered at the origin for each $i\ne j$, take $f_j$
to be the characteristic function of a ball of measure $\delta$ centered
at the origin, and let $\delta\to 0$.
Then $\tilde \Lambda(f_1,\cdots,f_m)$ has order of magnitude $\delta$,
while $\prod_i \|f_i\|_{L^{p_i}}$ has order of magnitude $\delta^{1/p_j}\ll
\delta$.

Valid inequalities can hold 
in Theorems~\ref{thm:discretegroups} and \ref{thm:amalgam}
with some exponents strictly less
than one, but they are always implied by stronger inequalities
already contained in those theorems.
More precisely,
if the inequality holds for some $m$-tuple $(p_1,\cdots,p_m)$,
then it also holds with each $p_i$ replaced by $\max(p_i,1)$. 
In the case of Theorem~\ref{thm:discretegroups}, that $p_j$ can be
replaced by $1$ if $p_j<1$ can be shown
by considering the case when the support of $f_i$ is a single point,
then exploiting linearity and symmetry.

\medskip
Two quite distinct investigations motivated our interest in these problems.
One derives from multilinear versions of the Kakeya-Nikodym maximal functions,
as will be explored in a forthcoming paper of the first, second, and fourth
authors.
A second motivator was work \cite{cltt} on multilinear operators with additional
oscillatory factors; see Proposition~\ref{prop:osc} and Corollary~\ref{cor:osc}
below. Further applications of Theorem~\ref{thm:global} to oscillatory integrals
will appear in a forthcoming paper \cite{christholmer}.

\section{An application to oscillatory integrals}

\begin{proposition} \label{prop:osc}
Let $m>1$.
For $1\le j\le m$
let $\ell_j:\reals^n\to\reals^{n_j}$ be surjective linear mappings. 
Let $P:\reals^n\to\reals$ be a polynomial.
Let $\varphi\in C^1_0(\reals^n)$ be a compactly supported,
continuously differentiable cutoff function.
For $\lambda\in\reals$ and $f_j\in L^{p_j}(\reals^{n_j})$
define
$\Lambda_{\lambda}(f_1,\cdots,f_m)
= \int_{\reals^n} e^{i\lambda P(x)}\prod_{j=1}^m f_j(\ell_j(x))\,\varphi(x)\,dx$. 
Suppose that
there exist 
$\delta>0$ and $C<\infty$
such that for all functions $f_j\in L^\infty$ and all $\lambda\in\reals$
\begin{equation}  \label{oscdecayforbounded}
|\Lambda_\lambda(f_1,\cdots,f_m)|\le C|\lambda|^{-\delta}
\prod_{j=1}^m \Norm{f_j}_{L^{\infty}}.
\end{equation}
Let
$(p_1,\cdots,p_m)\in[1,\infty]^m$, and suppose
that for every proper subspace $V\subset\reals^n$,
\begin{equation}
\codim_{\reals^n}(V) > \sum_j p_j\rp\codim_{\reals^{n_j}}(\ell_j(V)).
\end{equation}
Then there exist $\delta>0$ and $C<\infty$, depending on $(p_1,\cdots,p_m)$,
such that
\begin{equation} \label{oscillatorydecay}
|\Lambda_\lambda(f_1,\cdots,f_m)|\le C|\lambda|^{-\delta}
\prod_{j=1}^m \Norm{f_j}_{L^{p_j}}
\end{equation}
for all parameters $\lambda\in\reals$ and functions $f_j\in L^{p_j}(\reals^{n_j})$.
\end{proposition}

In the formulation of the hypothesis it is implicitly assumed
that the integral defining $\Lambda_\lambda(f_1,\cdots,f_m)$
converges absolutely for all functions $f_j\in L^{p_j}$;
thus by Theorem~\ref{thm:local} it is necessary that
$\codim_{\reals^n}(V)\le \sum_j p_j\rp\codim_{\reals^{n_j}}(\ell_j(V))$
for every subspace $V\subset\reals^n$.
The conclusion of Proposition~\ref{prop:osc} then
follows directly from Theorem~\ref{thm:local} by complex interpolation. 

A polynomial $P$ is said \cite{cltt}
to be nondegenerate, relative to the collection $\{\ell_j\}$
of mappings, if $P$ cannot be expressed as $P=\sum_j P_j\circ\ell_j$
for any collection of polynomials $P_j:\reals^{n_j}\to\reals$.

\begin{corollary} \label{cor:osc}
Let $\{\ell_j\},P,\varphi$ be as in Proposition~\ref{prop:osc}.
Suppose that $P$ is nondegenerate relative to $\{\ell_j\}$. 
Suppose that either
(i) $n_j=1$ for all $j$, $m<2n$, and the family $\{\ell_j\}$
of mappings is in general position,
or
(ii) $n_j=n-1$ for all $j$.
Let
$(p_1,\cdots,p_m)\in[1,\infty]^m$ 
and suppose that
for every proper subspace $V\subset\reals^n$,
$\codim_{\reals^n}(V) > \sum_j p_j\rp\codim_{\reals^{n_j}}(\ell_j(V))$.
Then there exists $\delta>0$ such that for any $\varphi\in C^1_0$
there exists $C<\infty$
such that for all functions $f_j\in L^{p_j}(\reals^{n_j})$,
\begin{equation*} 
|\Lambda_\lambda(f_1,\cdots,f_m)|\le C|\lambda|^{-\delta}
\prod_{j=1}^m \Norm{f_j}_{L^{p_j}}.
\end{equation*}
\end{corollary}
\noindent Here general position means that for any subset $S\subset\{1,2,\cdots,m\}$
of cardinality $|S|\le n$,
$\cap_{j\in S} \kernel(\ell_j)$ has dimension $n-|S|$.

\smallskip
By Theorems~2.1 and 2.2 of \cite{cltt}, the hypotheses imply
\eqref{oscdecayforbounded}.
Proposition~\ref{prop:osc} then implies the Corollary.

\section{Proof of sufficiency in Theorem~\ref{thm:global}}

We begin with the proof of sufficiency of the hypotheses
\eqref{homognecessary}, \eqref{globalmainhypothesis}
for the finiteness of the supremum in \eqref{maininequality}.
Necessity will be established in the next section.

The next definition is made for the purposes of the discussion
of Theorem~\ref{thm:global}; alternative notions of criticality
are appropriate for the other theorems.
\begin{definition}
A subspace $V\subset H$ is said to be critical if
\begin{equation}
\dim(V)=
\sum_j p^{-1}_j \dim(\ell_j(V)),
\end{equation}
to be supercritical if
the right-hand side is less than $\dim(V)$, 
and to be subcritical if the right-hand side is greater than $\dim(V)$.
\end{definition}
In this language,
the hypothesis \eqref{homognecessary} states that $V=H$ is critical,
while \eqref{globalmainhypothesis} states that no subspace of $H$ is supercritical.

\begin{proof}[Proof of sufficiency in Theorem~\ref{thm:global}]
The proof proceeds by induction on $\dim(H)$.
When $H$ has dimension one, necessarily $\dim(H_j)=1$ for all $j$. 
The hypothesis of the theorem in this case
is that $\sum_j p_j\rp=1$, and the conclusion 
is simply a restatement of H\"older's inequality for functions
in $L^{p_j}(\reals^1)$.

Suppose now that $\dim(H)>1$.
There are two cases. 
Case 1 arises when there exists some proper nonzero critical subspace 
$W\subset H$. 
The analysis then follows the pattern 
of \cite{finner} and \cite{CLL}.
Express $H = W^\perp\oplus W$ where $W^\perp$ is
the orthocomplement of $W$, with coordinates $y=(y',y'')\in W^\perp\oplus W$;
we will identify $(y',0)$ with $y'$ and $(0,y'')$ with $y''$.
Define 
$U_j\subset H_j$ to be 
\begin{equation}
U_j=\ell_j(W).
\end{equation}
Define $\tilde\ell_j=\ell_j|_W:W\to U_j$, which is surjective.
For $y'\in W^\perp$ and $x_j\in U_j$ define
\begin{equation}
g_{j,y'}(x_j) = f_j(x_j+\ell_j(y')).
\end{equation}
Then
\begin{equation}
f_j(\ell_j(y',y'')) 
= f_j(\ell_j(y') + \tilde\ell_j(y''))
= g_{j,y'}(\tilde\ell_j(y'')).
\end{equation}
Now
\begin{equation*}
\Lambda(f_1,\cdots,f_m)
=\int_{W^\perp}\int_W \prod_j f_j(\ell_j(y',y''))\,dy''\,dy'
\\
=\int_{W^\perp}\int_W \prod_j g_{j,y'}(\tilde\ell_j(y''))\,dy''\,dy',
\end{equation*}
so
\begin{equation}
\Lambda(f_1,\cdots,f_m)
= \int_{W^\perp}\tilde\Lambda(g_{1,y'},\cdots,g_{m,y'}) \,dy'
\end{equation}
where
\begin{equation}
\tilde\Lambda(g_1,\cdots,g_m) = \int_W \prod_j g_j(\tilde\ell_j(y''))\,dy''.
\end{equation}

We claim that 
\begin{equation}
\tilde\Lambda(g_1,\cdots,g_m) 
\le C\prod_j\|g_j\|_{p_j}.
\end{equation}
Since $W$ has dimension strictly less than  $\dim(H)$,
this follows from the induction hypothesis
provided that $W$ is critical and no subspace 
$V\subset W$ is 
supercritical, relative to the mappings $\tilde\ell_j$
and exponents $p_j$. But since $\tilde\ell_j$ is the restriction
of $\ell_j$ to $W$, 
this condition is simply the specialization of the original hypothesis from arbitrary
subspaces of $H$ to those subspaces contained in $W$,
together with the criticality of $W$ hypothesized in Case 1.
Thus
\begin{equation}
\Lambda(f_1,\cdots,f_m)
= \int_{W^\perp}\tilde\Lambda(g_{1,y'},\cdots,g_{m,y'}) \,dy'
\le C
\int_{W^\perp}\prod_j\|g_{j,y'}\|_{L^{p_j}(U_j)} \,dy'.
\end{equation}

We will next show how this last integral is another instance of the original problem,
with $H$ replaced by the lower-dimensional vector space
$W^\perp$.  For $z_j\in U_j^\perp$ define
\begin{equation}
F_j(z_j)=
\big(\int_{U_j} f_j(x_j+z_j)^{p_j}\,dx_j\big)^{1/p_j},
\end{equation}
recalling that $f_j\ge 0$,
with 
$F_j(z_j) = \esssup_{x_j\in U_j} f_j(x_j+z_j)$ if $p_j=\infty$.
Thus\footnote{If $U_j=\{0\}$ then the domain of $F_j$ is
$H_j$, and $F_j\equiv f_j$.
If $U_j=H_j$ then the domain of $F_j$ is $\{0\}$,
and $\|F_j\|_{p_j}$ is by definition $F_j(0) = \|f_j\|_{p_j}$.}
\begin{equation}
\|F_j\|_{L^{p_j}(U_j^\perp)}
=
\|f_j\|_{L^{p_j}(H_j)}.
\end{equation}

Denote by $\pi_{U_j^\perp}:H_j\to U_j^\perp$ 
and $\pi_{U_j}:H_j\to U_j$ the orthogonal projections.
Define
$L_j: W^\perp\to U_j^\perp$ by
\begin{equation}
L_j = \pi_{U_j^\perp}\circ\ell_j.
\end{equation} 
Decomposing
$\ell_j(y') = L_j(y')+u_j$ where $u_j=\pi_{U_j}(\ell_j(y'))$,
and making the change of variables $\tilde x_j = x_j+u_j$ in $U_j$,
gives
(if $p_j<\infty$)
\begin{multline}
\|g_{j,y'}\|_{L^{p_j}(U_j)}^{p_j}
=
\int_{U_j} |g_{j,y'}(x_j)|^{p_j}\,dx_j
=
\int_{U_j} |f_j(x_j+\ell_j(y'))|^{p_j}\,dx_j
\\
=
\int_{U_j} |f_j(x_j+u_j+L_j(y'))|^{p_j}\,dx_j
=
\int_{U_j} |f_j(\tilde x_j+L_j(y'))|^{p_j}\,d\tilde x_j
= 
F_j(L_j(y'))^{p_j}.
\end{multline}

Consequently we have shown thus far that
\begin{equation}
\Lambda(f_1,\cdots,f_m)
\le C
\int_{W^\perp}\prod_jF_j\circ L_j
\end{equation}
where 
$\|F_j\|_{L^{p_j}(U_j^\perp)} = \|f_j\|_{L^{p_j}(H_j)}$. 
Since $\ell_j:H\to H_j$ is surjective,
$H_j$ is spanned by $\ell_j(W)=U_j$ together with $\ell_j(W^\perp)$;
thus the orthogonal projection of $\ell_j(W^\perp)$ onto $U_j^\perp$ is all of
$U_j^\perp$; thus each $L_j:W^\perp\to U_j^\perp$ is surjective.

To complete the argument for Case 1 we need only show that
\begin{equation}
\int_{W^\perp}\prod_jF_j\circ L_j
\le C\prod_j\|F_j\|_{L^{p_j}(U_j^\perp)}.
\end{equation}
By induction on the ambient dimension, this follows from the next lemma,
which in the case when $\dim(H_j)=1$ for all $j$ appears in \cite{CLL}.
Although there are no additional complications in the general case, we include
a proof for the sake of completeness.

\begin{lemma}
Suppose that $H$ is critical, and has no supercritical subspaces.
Suppose that $W\subset H$ is a nonzero proper critical subspace. 
Define surjective linear transformations
$L_j=\pi_{\ell_j(W)^\perp}\circ \ell_j:W^\perp\to \ell_j(W)^\perp$. 
Then for any subspace $V\subset W^\perp$,
$\dim(V)\le\sum_j p_j\rp\dim(L_j(V))$.
\end{lemma}

\begin{proof}
Associate to $V$ the subspace $V+W\subset H$.
Since $V\subset W^\perp$,
$\dim(V+W)=\dim(V)+\dim(W)$.
Moreover,
for any $j$,
\begin{equation} \label{V+Wimagedimension}
\dim(\ell_j(V+W))
= \dim(L_j(V)) + 
\dim(\ell_j(W)),
\end{equation}
since 
$L_j = \pi_{\ell_j(W)^\perp}\circ\ell_j$.

Therefore 
\begin{align*}
\sum_j p_j^{-1}
\dim(L_j(V))
&=
\sum_j p_j^{-1}
\dim(\ell_j(V+W))
-
\sum_j p_j^{-1}
\dim(\ell_j(W))
\\
&=
\sum_j p_j^{-1}
\dim(\ell_j(V+W)) -\dim(W)
\\
& \ge \dim(V+W)-\dim(W)
=\dim(V),
\end{align*}
by the criticality of $W$ and subcriticality of $V+W$. 
Thus $V$ is not supercritical. 

When $V=W^\perp$, one has $V+W=H$,
whence $\sum_j p_j\rp \dim(\ell_j(V+W)) = \dim(V+W)$
since $H$ is assumed to be critical.
With this information the final inequality of the preceding display
becomes an equality, demonstrating that $W^\perp$ is critical.
\end{proof}

The proof of Case 1 of Theorem~\ref{thm:global} is complete. 
Turn next to Case 2, in which every nonzero proper subspace of $H$ is subcritical.
$\infty\rp$ is to be interpreted as zero throughout the discussion. 

Consider the set $K$ of all $m$-tuples $t=(t_1,\cdots,t_m)\in[0,1]^m$
such that relative to the exponents $p_j=t_j\rp$,
$H$ is critical and has no supercritical subspace.
Then $K$ equals the intersection of $[0,1]^m$ with a hyperplane and with various
closed half-spaces. Thus $K$ is convex and compact, 
whence it equals the closed convex hull of its extreme points. 

For any $t=(t_1,\cdots,t_m)\in [0,\infty)^m$, if \eqref{localhypothesis} holds, 
that is if
$\codim_H(V)\ge \sum_j t_j\codim_{H_j}(V)$ for all subspaces $V\subset H$,
then necessarily $t\in[0,1]^m$. 
Indeed, consider any index $i$ and let $V$ be the nullspace of $\ell_i$.
Then 
\begin{equation}
\dim(H_i) = \codim_H(V)
\ge \sum_j t_j\codim_{H_j}(\ell_j(V))
\ge t_i\codim_{H_i}\{0\}
= t_i\dim(H_i).\end{equation}

\eqref{localhypothesis} holds whenever $(t_j)=(p_j\rp)$
satisfies the hypotheses \eqref{homognecessary} and \eqref{globalmainhypothesis} 
of Theorem~\ref{thm:global}.
Consequently if $t$ is an extreme point of $K$, then some 
nonzero proper subspace of $H$ is critical relative to $t$,
or at least one coordinate $t_i$ equals $0$, 
or $m=1$ and $p_1=1$.
In the first subcase we are in Case 1, not Case 2. 
For the third subcase, see below.

In the second subcase, we may proceed by induction on the number $m$ of indices $j$,
for an inequality 
$\Lambda(f_1,\cdots,f_m)\le C\Norm{f_i}_{L^\infty}\prod_{j\ne i}\|f_j\|_{L^{p_j}}$
is equivalent to 
\begin{equation}
\Lambda(f_1,\cdots,f_{i-1},1,f_{i+1},\cdots,f_m)\le C\prod_{j\ne i}\|f_j\|_{L^{p_j}}.
\end{equation}
The hypotheses of Theorem~\ref{thm:global} are inherited by this
multilinear operator of one lower degree, acting on $\{f_j: j\ne i\}$,
whence the desired inequality follows by induction.

This induction is founded by the subcase where $m=1$, so that 
$\Lambda(f_1) = \int_H f_1\circ\ell_1$; moreover $p_1=1$.
Then $\ell_1:H\to H_1$ is surjective, so $\dim(H)\ge\dim(H_1)$.
The hypothesis $\dim(H)=p_1\rp\dim(H_1)\le \dim(H_1)$
thus forces $\ell_1:H\to H_1$ to be invertible, and $p_1\rp$ to equal $1$. 
Then $\Lambda(f_1) = c\int f_1$ for some finite constant $c$,
which is the desired result.
\end{proof}

\begin{remark}
When $\dim(H_j)=1$ for all $j$, every extreme point $(p_1\rp,\cdots,p_m\rp)$
of $K$ has each $p_j\rp\in\{0,1\}$ \cite{barthe},\cite{CLL}.
This is not the case in general;
in the Loomis-Whitney inequality for $\reals^n$, $K$ consists of a single point,
with $p_j =n-1$ for all $j$.
\end{remark}

\section{Proof of Theorem~\ref{thm:gut}} \label{gutsection}
Consider
$\int_{ \{y\in H: |\ell_0(y)|\le 1\} } \prod_{j=1}^m
f_j\circ\ell_j\,dy$ where the linear transformation $\ell_0$ has domain $H$
and range $H_0$ with $\dim(H_0)$ possibly equal to zero.
Thus some components of $y$ are constrained to a bounded set, while
the rest are free. 
Set 
\begin{equation}\scriptv=\kernel(\ell_0);\end{equation}
the component of $y$
lying in $\scriptv$ is completely unconstrained, 
while the component in $\scriptv^\perp$ is constrained to a bounded set.

\begin{proof}[Proof of necessity of \eqref{firstnec} and \eqref{secondnec}]
For any subspace $V\subset H$
define $V_\bg= V \cap\scriptv$
and $V_\smll= V\ominus V_\bg$, so that 
$V = V_\smll\oplus V_\bg$.
Let $r\le 1\le R$ be arbitrary. 
Define $f_j=f_j(x_j)$ to be the characteristic function of 
the region $S_j$ where $|x_j|\le R$ if $x_j\in\ell_j(V_\bg)$,
$|x_j|\le 1$ if $x_j\in \ell_j(V)\cap (\ell_j(V_\bg))^\perp$,
and $|x_j|\le r$ if $x_j\in (\ell_j(V))^\perp$.

Let $c_0>0$ be a small constant, independent of $r,R$, and
define $S\subset H$ to be the set of all $y$
such that $|y|\le c_0 r$ if $y\in V^\perp$,
$|y|\le c_0$ if $y\in V_\smll$,
and $|y|\le c_0 R$ if $y\in V_\bg$.
Then provided $c_0$ is chosen sufficiently small,
$y\in S\Rightarrow f_j(\ell_j(y))=1$ for all indices $j$.
Indeed, if $y\in V^\perp$ then $|\ell_j(y)|\le C |y|\le C c_0 r$,
so $\ell_j(y)\in S_j$.
If $y\in V_\smll$ then $|\ell_j(y)|\le C |y|\le C c_0$,
so since $\ell_j(y)\in\ell_j(V)$, $\ell_j(y)\in S_j$.
Finally if $y\in V_\bg$ then $|\ell_j(y)|\le C c_0 R$,
which implies that $\ell_j(y)\in S_j$ since $\ell_j(y)\in \ell_j(V_\bg)$.

Moreover $y\in S\Rightarrow |\ell_0(y)|\le 1$.
Therefore
\begin{equation}
\tilde\Lambda_{\text{loc}}(\{f_j\})
\ge |S| 
\sim 
R^{\dim(V_\bg)}
\cdot
r^{\codim_H(V)}
\end{equation}
while
\begin{equation}
\|f_j\|_{p_j}
\sim
R^{p_j\rp \dim(\ell_j(V_\bg))}
r^{p_j\rp \codim_{H_j}(\ell_j(V))}.
\end{equation}
Suppose that the ratio 
$\tilde\Lambda_{\text{loc}}\big(\{f_j\}\big) / \prod_j\|f_j\|_{p_j}$
is bounded uniformly as a function of $r,R$.
By letting $r\to 0$, 
we conclude that 
$\dim(V_\bg)\le\sum_j p_j\rp\dim(\ell_j(V_\bg)$.
Letting $R\to\infty$ gives 
$\codim_H(V)\ge \sum_j p_j\rp \codim_{H_j}(\ell_j(V))$.
\end{proof}

The following lemma will be used in the proof of Theorem~\ref{thm:gut}.
\begin{lemma} \label{localthmlemma}
Suppose that 
$\codim_H(V)\ge \sum_j p_j\rp\codim_{H_j}(\ell_j(V))$ 
for every
subspace $V\subset H$, and that $W\subset H$ is a subspace
satisfying $\codim_H(W) = \sum_j p_j\rp\codim_{H_j}(\ell_j(W))$.
Then  for any subspace $V\subset W$,
$\codim_W(V)\ge \sum_j p_j\rp\codim_{\ell_j(W)}(\ell_j(V))$.
Likewise for any subspace $V\subset W^\perp$, 
$\codim_{W^\perp}(V)\ge \sum_j p_j\rp\codim_{\ell_j(W)^\perp}(L_j(V))$.
\end{lemma}

\begin{proof}
For the first conclusion,
\begin{multline}
\codim_W(V) 
= 
\dim(W)-\dim(V) 
= \codim_H(V) - \codim_H(W)
\\
\ge  \sum_j p_j\rp \codim_{H_j}(\ell_j(V)) - \sum_j p_j\rp \codim_{H_j}(\ell_j(W))
\\
= \sum_j p_j\rp (\dim(\ell_j(W))-\dim(\ell_j(V)))
= \sum_j p_j\rp \codim_{\ell_j(W)}(\ell_j(V)).
\end{multline}

For the second conclusion,
\begin{equation} 
\begin{split}
\codim_{W^\perp}(V)
&= \dim(H)-\dim(W)-\dim(V)
\\
& = \codim_H(V+W)
\\
&\ge \sum_j p_j\rp \codim_{H_j}(V+W)
\\
&= \sum_j p_j\rp \big(\dim(H_j)-\dim(\ell_j(W))-\dim(L_j(V))\big)
\\
&= \sum_j p_j\rp \big(\dim(L_j(W^\perp))-\dim(L_j(V))\big).
\\
&= \sum_j p_j\rp \codim_{L_j(W^\perp)}(L_j(V)).
\end{split}
\end{equation}
The identity $\dim(H_j) = \dim(\ell_j(W))
+\dim(L_j(W^\perp))$ used to obtain the final line is
\eqref{V+Wimagedimension} specialized to $V=W^\perp$.
\end{proof}

\begin{proof}[Proof of sufficiency in Theorem~\ref{thm:gut}]
The proof follows the inductive scheme of the proof 
of Theorem~\ref{thm:global}. 
To simplify notation set $t_j=p_j\rp\in[0,1]$.
Case 1 now breaks down into two subcases.
Case 1A arises when there exists a nonzero proper subspace 
$W$ of $H$ that is contained in $\scriptv$ and 
is critical in the sense of \eqref{firstnec}, that is,\footnote{
All summations with respect to $j$ are taken over $1\le j\le m$.}
$\sum_j t_j\dim(\ell_j(W)) = \dim(W)$.

With coordinates $(y',y'')$ for $W^\perp\oplus W$,
$\ell_0$ is independent of $y''$, and for every subspace $V\subset W$,
$\sum_j t_j\dim(\ell_j(V))\ge\dim(V)$ by \eqref{firstnec}. Thus the collection
of mappings $\{\ell_j|_W \}$ satisfies the hypothesis of 
Theorem~\ref{thm:global}, 
whence
$\int_W \prod_j f_j\circ\ell_j(y',y'')\,dy''
\le C \prod_j F_j(y')$ 
where $\|F_j\|_{L^{p_j}(W^\perp)}\le C\|f_j\|_{L^{p_j}(H_j)}$.

It remains to bound $\int_{W^\perp} \chi_B\circ\ell_0(y',0)
\prod_j F_j\circ L_j(y')\,dy'$, where $B$ denotes the characteristic
function of a ball of finite radius.
Theorem~\ref{thm:gut} can be invoked by induction on the ambient dimension,
provided that \eqref{firstnec} and \eqref{secondnec} hold for the 
data $W^\perp,\scriptv\cap W^\perp,\{U_j^\perp,L_j,p_j\}$.
We will write \eqref{firstnec}$_{H}$, \eqref{firstnec}$_{W}$,
and \eqref{firstnec}$_{W^\perp}$ to distinguish between this
hypothesis for the three different data that arise in the discussion;
likewise for \eqref{secondnec}.

\eqref{secondnec}$_W$ is the condition that
$\codim_{W^\perp}(V)
\ge \sum_j t_j \codim_{L_j(W^\perp)}(L_j(V))$
for every subspace $V\subset W^\perp$,
which is the second conclusion of Lemma~\ref{localthmlemma}.
\eqref{firstnec}$_W$ is the condition
\begin{equation} \label{firstnecalt1}
\dim(V)\le \sum_j t_j\dim(L_j(V)) 
\ \text{ for all subspaces $V\subset \scriptv\cap W^\perp$}.
\end{equation}
Since $V,W$ are both contained in $\scriptv$ so is $V+W$,
so $\sum_j t_j \dim(\ell_j(V+W))\ge \dim(V+W) = \dim(V)+\dim(W)$
by \eqref{firstnec}$_H$.
This together with the previously established identity
$\dim(\ell_j(V+W)) = \dim(\ell_j(W)) + \dim(L_j(V))$
and the criticality condition $\sum_j t_j\dim(\ell_j(W)) = \dim(W)$
yields \eqref{firstnecalt1}.
Thus Case 1A is treated by applying Theorem~\ref{thm:global} for $W$
and the induction hypothesis for $W^\perp$.

Case 1B arises when there exists a nonzero proper subspace $W\subset H$
that is critical in the sense of \eqref{secondnec}, that is,
$\codim_H(W) = \sum_j t_j \codim_{H_j}(\ell_j(W))$. The analysis follows the same 
inductive scheme.
Lemma~\ref{localthmlemma} guarantees that \eqref{secondnec}$_W$ holds, while
\eqref{firstnec}$_W$  is simply the specialization
of \eqref{firstnec}$_H$ to subspaces $V\subset W\cap \scriptv$.
Thus Theorem~\ref{thm:gut} may be applied by induction to
$W,\{\ell_j(W),\ell_j|_W,p_j\}$.

This reduces matters
to $\int_{W^\perp\cap \{|L_0(y')|\le 1\}} \prod_j F_j\circ L_j\,dy'$,
where the nullspace $\tilde V$ of $L_0$
is the set of all $y'\in W^\perp$ for which there exists $y''\in W$
such that $\ell_0(y',y'')=0$;
thus the subspace $\scriptv\subset H$ 
is now replaced by $\pi_{W^\perp}\scriptv\subset W^\perp$.

Now it is natural to expect to use \eqref{firstnec}$_H$
to establish \eqref{firstnec}$_{W^\perp}$,
but the latter pertains to certain subspaces not contained in 
$\scriptv$, about which the former says nothing.
Luckily the inequality in \eqref{firstnecalt1} holds for arbitrary subspaces
$V\subset W^\perp$, not merely those contained 
in $\pi_{W^\perp}\scriptv$.
Indeed,
\begin{align*}
\sum_j t_j\dim(L_j(V))
&=\sum_j t_j\dim(\ell_j(V+W)) -\sum_j t_j\dim(\ell_j(W))
\\
&= \sum_j t_j \codim_{H_j}(\ell_j(W))
- \sum_j t_j \codim_{H_j}(\ell_j(V+W))
\\
&= \codim_H(W)
- \sum_j t_j \codim_{H_j}(\ell_j(V+W))
\\
&\ge \codim_H(W)
- \codim_H(V+W)
\\
&=\dim(V).
\end{align*}

The assumption that $W$ is critical in the
sense that equality holds in \eqref{secondnec}$_H$ implies 
\eqref{secondnec}$_W^\perp$, 
by the second conclusion of Lemma~\ref{localthmlemma}.
Thus by induction on the dimension, Theorem~\ref{thm:gut} may be applied
to the integral over $W^\perp$, concluding the proof for Case 1B. 

Case 2 arises when no subspace $W$ is critical in either sense.
Consider the set $K\subset[0,1]^m$ of all $(t_1,\cdots,t_m)$ such that
$\sum_j t_j \dim(\ell_j(V))\ge\dim(V)$ for all subspaces $V\subset\scriptv
=\kernel(\ell_0)$,
and
$\codim_H(V)\ge \sum_j t_j\codim_{H_j}(\ell_j(V))$ for all subspaces
$V \subset H$.
It suffices to prove that $\int_{H} \chi_B\circ\ell_0\prod_{j\ge 1}
f_j\circ\ell_j \le C\prod_j\|f_j\|_{q_j}$ for every extreme point
$(t_1,\cdots,t_m)$ of $K$, where $q_j=t_j\rp$.
Consider such an extreme point.
If there exists a nonzero proper subspace $V\subset\scriptv$ that is critical in the
sense that $\sum_j t_j \dim(\ell_j(V))=\dim(V)$,
or a nonzero proper subspace $V\subset H$ that is critical in the sense that 
$\codim_H(V)= \sum_j t_j\codim_{H_j}(\ell_j(V))$,
then Case 1A or Case 1B apply.

There are other cases in which equality might hold in 
\eqref{firstnec} or \eqref{secondnec}, besides those subsumed under Case 1.
If equality holds for $V=\{0\}$ in \eqref{secondnec} with $p_j\rp=t_j$,
then $\dim(H)=\sum_j t_j \dim(H_j)$, 
which is the first hypothesis of Theorem~\ref{thm:global}.
In conjunction with \eqref{secondnec} this
implies that \eqref{firstnec} holds for every subspace $V\subset H$,
which is the second hypothesis of Theorem~\ref{thm:global}. Therefore 
the conclusion \eqref{gutinequality} of Theorem~\ref{thm:gut}
holds without the restriction $|\ell_0(y)|\le 1$ in the integral, 
by Theorem~\ref{thm:global}.

If on the other hand $H=\scriptv=\kernel(\ell_0)$ 
and equality holds for $V=H$ in \eqref{firstnec} with $p_j\rp=t_j$,
then $\dim(H)=\sum_j t_j \dim(H_j)$, so 
Theorem~\ref{thm:global} applies once more.

Therefore matters reduce to the
case where equality holds in \eqref{firstnec} for  no subspace of $\scriptv$
except $V=\{0\}$, and where furthermore equality holds in \eqref{secondnec}
for no subspace of $H$ except for $V=H$ itself.
Equality always holds in both of those cases, so they play no part
in defining $K$.

$t$ satisfies
$\codim_H(V)\ge \sum_j t_j\codim_{H_j}(V))$ for every subspace $V\subset H$.
Therefore as in Case 2 of the proof of Theorem~\ref{thm:global},
every remaining extreme point
$(t_1,\cdots,t_m)$ of $K$ must have $t_i=0$ for at least one index $i$.

By induction on $m$, it therefore suffices to treat the case $m=1$, with $p_1=\infty$.
By \eqref{firstnec} applied to $V = \kernel(\ell_0)$, $\dim(\kernel(\ell_0))\le 0\dim(H_1)=0$,
so $\ell_1$ has no kernel. Therefore the restriction $|\ell_0(y)|\le 1$
constrains $y$ to a bounded region, whence $\int_{|\ell_0(y)|\le 1} f_1\circ\ell_1(y)\,dy
\le C \|f_1\|_{L^\infty}$ for some finite constant $C$.
\end{proof}

\section{Proof of Theorem~\ref{thm:discretegroups}}
This proof contains no new elements, so will merely be outlined.
We denote the identity element of a group by $0$.
Recall that
if $H_1,H_2$ are subgroups of a finitely generated
discrete Abelian group $G$, and if $H_1\cap H_2=\{0\}$, 
then $\rank(H_1+H_2)=\rank(H_1)+\rank(H_2)$.
Likewise if $H'$ is a subgroup of the quotient group $G/H$ then
$\rank(H)+\rank(H')$ equals $\rank(\pi\rp(H'))$ where $\pi:G\to G/H$
is the natural projection. A finitely generated Abelian group is finite
if and only if its rank is zero.

Let groups $G,G_j$, homomorphisms $\varphi_j$, and exponents $p_j$ satisfy the 
hypotheses of Theorem~\ref{thm:discretegroups}. Consider first the case where
there exists a subgroup $G'\subset G$, satisfying $0<\rank(G')<\rank(G)$,
that is critical in 
the sense that $\sum_j p_j\rp \rank(\varphi_j(G')) = \rank(G')$.
Define $G'_j = \varphi_j(G')\subset G_j$.
Since every subgroup of $G$ inherits the hypothesis of the theorem,
we may conclude by induction on the rank that
\begin{equation}
\sum_{y\in G'} \prod_j f_j\circ\varphi_j(y)\le C\prod_j \|f_j\|_{\ell^{p_j}(G'_j)}.
\end{equation}
Define 
$F_j\in \ell^{p_j}(G_j/G'_j)$ by 
\[F_j(x+G'_j) = (\sum_{z\in G'_j}|f_j(x+z)|^{p_j})^{1/p_j}.\] 
Then $\|F_j\|_{\ell^{p_j}(G_j/G'_j)}\le C \|f_j\|_{\ell^{p_j}(G_j)}$.
Define homomorphisms $\psi_j: G/G'\to G_j/G'_j$
by composing $\varphi_j$ with the quotient map from $G_j$ to $G'_j$.
Then 
\begin{equation}
\sum_{y\in G} \prod_j f_j\circ\varphi_j(y)
= \sum_{x\in G/G'} \sum_{z\in G'} f_j\circ\varphi_j(x+z)
\le \sum_{x\in G/G'} F_j\circ\psi_j(x).
\end{equation}
It suffices to show that the homomorphisms $\psi_j$ inherit the
hypothesis of Theorem~\ref{thm:discretegroups}, which may then
be applied by induction on the rank to yield the desired bound
$O(\prod_j\|F_j\|_{\ell^{p_j}})$.
This hypothesis is verified using the criticality of $G'$ and
the additivity of ranks, just as in the proof of Theorem~\ref{thm:global}.

There remains the case in which no critical subgroup $G'$ of strictly smaller
but strictly positive rank exists. Once again we consider the compact convex
set $K$ of all $(q_1\rp,\cdots,q_N\rp)\in[0,1]^m$ 
for which $\rank(H)\le \sum_j q_j\rp \rank(\varphi_j(H))$
for all subgroups $H\subset G$, and it suffices to
prove that $\sum_{y\in G} \prod_j f_j\circ\varphi_j(y) \le C\prod_j \|f_j\|_{q_j}$
for all extreme points $(q_1\rp,\cdots,q_N\rp)$ of $K$.

If $(q_1\rp,\cdots,q_N\rp)$ is an extreme point 
then either 
$\sum_j q_j\rp \rank(\varphi_j(G')) =\rank(G')$ for some subgroup $G'$ 
satisfying $0<\rank(G')<\rank(G)$, or $q_j\rp\in\{0,1\}$
for all indices $j$, or $\rank(G)= \sum_j q_j\rp\rank(\varphi_j(G))$ 
and $q_j\rp\in\{0,1\}$ for all but at most one index $j$. 
In the first case we are in the critical case treated above. 

Suppose 
that $(q_1\rp, \cdots,q_N\rp)\in K$ and $q_j\in\{0,1\}$ for all $j$.
Let $S=\{j: q_j\rp=1\}$, and consider the subgroup 
$G'=\cap_{j: q_j=1} \kernel(\varphi_j)$.
The hypothesis \eqref{discretehypothesis} 
states that $0 = \sum_{j\in S} \rank(\varphi_j(G')) \ge \rank(G')$, so $G'$ is finite.
If $f_j$ is the characteristic function of a single point $z_j$ for each $j\in S$,
then
$\sum_{y\in G} \prod_{j\in S} f_j\circ\phi_j(y)$ equals the cardinality
of $\{y: \phi_j(y)=z_j\ \forall j\in S\}$, which is $\le |G'|$.
The inequality then follows for arbitrary functions by multilinearity.

Suppose finally that $q_i\rp\in (0,1)$, $q_j\rp=1$ if and only if $j\in S$,
and $q_j\rp=0$ if neither $j\in S$ nor $j=i$.
Let $S=\{j: q_j\rp=1\}$ and
consider $G' = \cap_{j\in S} \kernel(\varphi_j)$. The hypothesis
\eqref{discretehypothesis} states that
$\rank(G')\le q_i\rp\rank(\varphi_i(G')) + \sum_{j\in S}\rank(\varphi_j(G'))
= q_i\rp \rank(\varphi_i(G'))$; the right-hand side is necessarily
$\le q_i\rp \rank(G')$, which is strictly less than $\rank(G')$
unless $\rank(G')=0$; hence $\rank(G')$ must vanish.
Therefore for any nonnegative functions, 
\begin{equation*}
\sum_{y\in G} \prod_{j} f_j\circ\varphi_j(y)
\le C\prod_{j\in S}\|f_j\|_{\ell^1}\prod_{j\notin S} \|f_j\|_{\ell^\infty},
\end{equation*}
as in the preceding paragraph.
Since $\|f_i\|_\infty\le \|f_i\|_{q_i}$, this completes the proof.
\qed

The proof of the variant Theorem~\ref{thm:amalgam} is nearly identical
to that of Theorem~\ref{thm:discretegroups} and is left to the reader.
Likewise the proofs of the necessity of the hypotheses in both theorems,
which are simplifications of the reasoning shown above for their
continuum analogues, are omitted.

\section{Variants based on product structure}  \label{section:finner}

A variant of our results, based on combinatorial rather than
linear algebraic or group theoretic structure, 
has been obtained earlier by Finner \cite{finner}; see also 
\cite{friedgut} for a discussion of some special cases from another
point of view.
Let $\{(X_i,\mu_i)_{i \in I}\}$
be a finite collection of measure spaces,
and let $(X,\mu)=\prod_{i\in I}(X_i,\mu_i)$ be their product. 
Let $J$ be another finite index set.
For each $j\in J$, let $S_j$ be some nonempty subset of $I$.
Let $Y_j = \prod_{i\in S_j} X_i$, equipped with the associated product
measure, and let $\pi_j:X\to Y_j$ be the natural projection map.
Let $f_j: Y_j\to[0,\infty]$ be measurable.
To avoid trivialities, we assume throughout the discussion that $I,J$ are nonempty
and that $\mu(X)$ is strictly positive.
Define 
\begin{equation}
\Lambda(f_j)_{j\in J}
= \int_X \prod_{j\in J} f_j\circ\pi_j\,d\mu. 
\end{equation}
Denote by $|\cdot|$ the cardinality of a finite set.

Let $p_j\in [1,\infty]$ for each $j\in J$.
Finner's theorem then asserts that if
\begin{equation} \label{finnerhypothesis}
1=\sum_{j: i\in S_j}p_j\rp \text{ for all $i\in I$}
\end{equation}
then
\begin{equation} \label{finnerconclusion}
\Lambda(f_j)_{j\in J}
\le
\prod_{j\in J}\Norm{f_j}_{L^{p_j}(Y_j)}.
\end{equation}
The hypothesis \eqref{finnerhypothesis}
can be equivalently restated as
\begin{equation} \label{finnercriticality}
|K| = \sum_{j\in J} p_j\rp |S_j\cap K| 
\text{ for every subset $K\subset I$,}
\end{equation}
or again as the conjunction of
$|I| = \sum_{j\in J} p_j\rp |S_j|$
and 
$|K| \le \sum_{j\in J} p_j\rp |S_j\cap K|$ for every $K\subset I$.
The analogue of a subspace is now a subset $K\subset I$,
and the analogue of criticality is \eqref{finnercriticality};
the inequality need not hold, in general, unless every subset  is critical. 
This contrasts with the situation treated by Carlen, Lieb, and Loss \cite{CLL}
and in Theorem~\ref{thm:global},
where generic subspaces will be subcritical even if critical subspaces exist.

When each space $X_i$ is some Euclidean space
equipped with Lebesgue measure, the hypotheses in this last form are precisely those
of Theorem~\ref{thm:global}, specialized to this limited class of 
linear mappings. 
A special case is the Loomis-Whitney inequality
\begin{equation*}
\int_{\reals^n} \prod_{j=1}^n f_j\circ\pi_j(x)\,dx
\le \prod_{j=1}^n \Norm{f_j}_{L^{n-1}},
\end{equation*}
where $\pi_j:\reals^n\to\reals^{n-1}$ is the mapping that forgets
the $j$-th coordinate.

Our next result is analogous to a unification of Theorems~\ref{thm:gut}
and \ref{thm:amalgam}.
We say that a measure space $(X,\mu)$ is atomic if there exists
$\delta>0$ such that $\mu(E)\ge\delta$ for every
measurable set $E$ having strictly positive measure.
\begin{proposition} \label{prop:finnergut}
Suppose that the index set $I$ is a disjoint union $I = I_0\cup I_\infty\cup I_\star$,
where $X_i$ is a finite measure space for each $i\in I_0$,
is atomic for each $i\in I_\infty$, 
and is an arbitrary measure space for each $i\in I_\star$.
Then a sufficient condition for the inequality \eqref{finnerconclusion} is that
\begin{alignat}{2}
&1\ge &&\sum_{j: i\in S_j}p_j\rp \text{ for all $i\in I_0$}
\label{finnerlocalhypothesis}
\\
&1\le &&\sum_{j: i\in S_j}p_j\rp \text{ for all $i\in I_\infty$}
\label{finneratomichypothesis}
\\
&1 =  &&\sum_{j: i\in S_j}p_j\rp \text{ for all $i\in I_\star$}.
\label{finnerhypothesisagain}
\end{alignat}
\end{proposition}
\noindent That these sufficient conditions are also necessary, in general,
is a consequence of the necessity of the hypotheses of Theorem~\ref{thm:gut}.

\begin{remark} \label{finnerremark}
Consider the case where
each $X_i$ is a finite measure space. 
If $(p_j)_{j\in J}$ satisfies the hypothesis \eqref{finnerhypothesis},
and if $q_j\ge p_j$ for all $j\in J$,
then 
$\Lambda(f_j)_{j\in J} 
\le C\prod_j\Norm{f_j}_{p_j}
\le C'\prod_j\Norm{f_j}_{q_j}$
by Finner's theorem and H\"older's inequality.
However, there are situations\footnote{The special case of
Proposition~\ref{prop:finnergut} in which all $X_i$ are finite measure spaces
is stated in \cite{finner}, p.\ 1898, but no proof is given.}
in which $(q_j)_{j\in J}$ satisfies
\eqref{finnerlocalhypothesis} 
yet there exists no $(p_j)_{j\in J}$ satisfying \eqref{finnerhypothesis}
with $q_j\ge p_j$ for all $j\in J$.

To construct an example, begin with any situation where there is an
extreme point $(q_j\rp)_{j\in J}$ of $K=\{(t_j)_{j\in J}\in[0,1]^J:
1 = \sum_{j: i\in S_j} t_j \text{ for all } i\in I\}$,
such that $q_j\rp<1$ for all $j$; for instance, the Loomis-Whitney example.
Augment $I$ by adding a single new index $i'$,
choose one index $j'$ already in $J$, and replace $S_{j'}$ by $S_j\cup \{i'\}$,
while keeping $S_j$ unchanged for all $j\ne j'$.
Thus $\sum_{j: i'\in S_j} q_j\rp = q_{j'}\rp<1$; 
$(q_j)_{j\in S}$ satisfies \eqref{finnerlocalhypothesis}. 
However no $(p_j)_{j\in J}$. For if $p_j\ge q_j\rp$ for all $j$ with
strict inquality for some index $k$, choose some $i\in S_k$.
Then $\sum_{j: i\in S_j} p_j\rp >\sum{j: i\in S_j}q_j\rp =1$,
so that \eqref{finnerlocalhypothesis} fails for $(p_j)_{j\in S}$.
\end{remark}

Proposition~\ref{prop:finnergut} can be proved by repeating Case 1 of the proofs of
Theorems~\ref{thm:global} and \ref{thm:gut}, arguing by induction on $|I|$,
and integrating with respect to the $m$-th coordinate in $\prod_{i\in I}X_i$
while all other coordinates are held constant. 
The basis case $m=1$ is H\"older's inequality.
Indeed, this is the argument given in \cite{finner} for the special case when $I=I_\star$.

Alternatively, when $I_0$ is empty,\footnote{To treat the general case in this way
would require a unification of Theorems~\ref{thm:gut} and \ref{thm:amalgam}
analogous to Proposition~\ref{prop:finnergut}. We see no obstruction to such
a result.} 
Proposition~\ref{prop:finnergut}
can be reduced to the case where each $X_i$ is $\reals^1$ equipped with Lebesgue measure, 
by approximating general functions by finite linear combinations of characteristic
functions of product sets, and then embedding any particular situation measure-theoretically
into a (product of copies of) $\reals^1$. The inequality \eqref{finnerconclusion} 
then follows from an application of Theorem~\ref{thm:global}.

\end{document}\end